\documentclass{article}
\usepackage{amssymb}
\usepackage{amsmath}
\usepackage{epsfig,psfrag,verbatim,amsfonts}
\usepackage[tight]{subfigure}
\usepackage{epstopdf}
\begin{document}

\newtheorem{theorem}{Theorem}
\newtheorem{lemma}{Lemma}
\newtheorem{corollary}{Corollary}
\newtheorem{proposition}{Proposition}
\def\squarebox#1{\hbox to #1{\hfill\vbox to #1{\vfill}}}
\newcommand{\qed}{\hspace*{\fill}
\vbox{\hrule\hbox{\vrule\squarebox{.667em}\vrule}\hrule}\smallskip}
\newcommand{\T}{{\mathcal{T}}}
\newcommand{\ZZ}{{\mathbb{Z}}}
\newcommand{\reals}{{\mathbb{R}}}
\renewcommand{\P}{{\mathbf{P}}}
\newcommand{\eps}{{\varepsilon}}
\renewcommand{\v}{{\tt v}}
\newcommand{\BBB}{\mbox{${\mathbb{B}}$}}
\title{{\bf
Localization for controlled random walks
and martingales}}

\author{Ori Gurel-Gurevich\thanks{University of British Columbia}\and
Yuval Peres\thanks{Microsoft Research}
\and
Ofer Zeitouni\thanks{Weizmann Institute of Science and Courant institute.
Research partially supported by
NSF grant \#DMS-1106627 and  by a
grant from the Israel Science Foundation.}}
\date{September 16, 2013}
\maketitle
\begin{abstract}
 We consider controlled random walks that are martingales with uniformly bounded increments and nontrivial
  jump probabilities and show that
  such walks can be constructed so that $P(S_n^u=0)$ decays at polynomial
 rate $n^{-\alpha}$ where $\alpha>0$ can be arbitrarily small. We also
 show, by means of a general delocalization lemma for martingales, which
 is of independent interest,
 that slower than polynomial decay is not possible.
\end{abstract}

%AMS Subject classification: primary 60K37, 60F05. Secondary 60J80, 82C41.

\section{Introduction and statement of results}
\label{sec-introduction}
\setcounter{figure}{0}
Consider a discrete time martingale $\{M_i\}_{i\geq 0}$
adapted to a filtration ${\cal F}_i$
whose increments are uniformly bounded by $1$, i.e. $|M_{i+1}-M_i|\leq 1$,
and such that $P(|M_{i+1}-M_i|=1\mid {\cal F}_i)>c>0$.
%its quadratic variation
%process satisfies the bounds $c_1 n \langle M\rangle_n\leq c_2 n$ for some
%positive constants $c_1,c_2>0$.
It is folklore that in many respects, such a martingale should be well
approximated by Brownian motion. In particular, one would expect that
$P(|M_n|\leq 1)$ should be of order $n^{-1/2}$.

Our goal in this paper is to point out that this naive expectation is completely
wrong.
We will frame this in the language of controlled processes below, but a
corollary of our main result, Theorem \ref{theo-main} below, is the following.
\begin{corollary}
  \label{cor-main}
  For any $\alpha>0$ there exist  $\beta>0$ and $c>0$ so that for any $n>0$
  there exists
  an ${\cal F}_i$-adapted
  discrete time martingale $\{M_i\}_{i\geq 0}$
  with
  $|M_{i+1}-M_i|\leq 1$ and $P(|M_{i+1}-M_i|=1\mid {\cal F}_i)>\beta$
  such that
  $$ P(|M_n|\leq 1)\geq c n^{-\alpha} \,.$$
\end{corollary}

\begin{figure}
  \epsfig{file=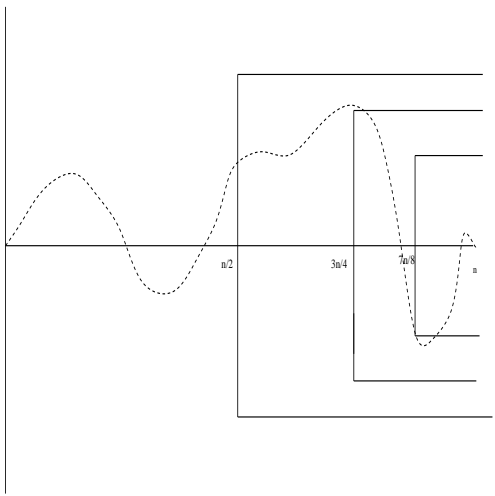,height=6.5cm, width=6.5cm}
\begin{centering}
\caption{The sets $D_i$ and their crossings by a trajectory with $M_n=0$.}
\end{centering}\label{fig-1}
\end{figure}

Corollary \ref{cor-main} can be viewed as a localization lemma. A
complementary delocalization estimate was obtained by de la Rue \cite{TDLR}.
We provide a
different proof to a strengthened version of his results.
\begin{theorem}
\label{theo-gen}
For any $\delta\in (0,1]$ and $\beta\in [0,1/2)$ there exist
$C=C(\delta,\beta)<\infty$ and $\alpha=\alpha(\delta,\beta)>0$ so that
the following holds.

If $\{M_i\}_{i\geq 0}$ is
a discrete time martingale (with respect to a filtration ${\cal F}_i$)
satisfying $E( (M_{i+1}-M_i)^2|{\cal F}_i)\in [\delta,1]$ and
$|M_{i+1}-M_i|\leq n^\beta$ a.s., then
\begin{equation}
\label{gen-UB}
%\limsup_{n\to\infty}  n^{\alpha(q)}
\sup_z P(|M_n-z|\leq n^\beta )<C n^{-\alpha}\,.
  \end{equation}
  \end{theorem}
  The heart of the 
  proof of Theorem \ref{theo-gen} uses a sequence of entrance times to a
  space-time region, which may be of independent interest
  (see Figure \ref{fig-1} for a graphical depiction).

Our interest in these questions arose while one of us was working on \cite{DLP}.
Charlie Smart then kindly pointed out  \cite{Smart}
that the continuous
time results in \cite{Atrok} and \cite{BaS} concerning the viscosity
solution of certain optimal control problems could be adapted to the discrete
time setting (using \cite{BS}) in order to show an integrated version of
Corollary \ref{cor-main}, namely that for any fixed
$\beta,\gamma>0$ a
martingale $\{M_i\}$ as in the lemma could be constructed so that
for all $\delta$ small,
\begin{equation}
  \label{eq-222}
  P(|M_n|\leq \delta \sqrt{n})\geq \gamma \delta \,.
\end{equation}
(Note that $\gamma$ can be taken arbitrarily large, for $\beta$ fixed.
The estimate
\eqref{eq-222}
is in contrast with the expected linear-in-$\delta$
behavior
one might naively expect from diffusive scaling.)
This then raises the question
of whether a local version of this result could be obtained,
and our goal in this short note is to answer that in the affirmative.

We phrase some of
our results in the language of \textit{controlled random walks}.
Fix a parameter $q\in [0,1)$.
Consider a controlled simple random walk $\{S_i^u\}_{i\geq 0}$,
defined as follows. Let $S_0=0$ and
let ${\cal F}_i=\sigma(S_0,S_1,\ldots,S_i)$ denote the sigma-field generated
by the process up to time $i$. A \textit{$q$-admissible control}
is a collection of random variables $\{u_i\}_{i\geq 0}$ satisfying the following
conditions:
\begin{enumerate}
  \item
a) $u_i\in [0,q]$, a.s..
\item
b) $u_i$ is ${\cal F}_i$-adapted.
\end{enumerate}

Let ${\cal U}_q$ denote the set of all $q$-admissible controls.
For $u\in {\cal U}_q$,
the controlled simple random walk $\{S_i^u\}_{i\geq 0}$ is determined by
the equation
\begin{equation}
  P(S_{i+1}^u=S_i^u +\Delta| {\cal F}_i)=
  \left\{\begin{array}{ll}
    u_i, & \Delta=0\\
    (1-u_i)/2,& \Delta=\pm 1\,.
  \end{array}
  \right.
\end{equation}
Of course, $\{S_i^u\}_{i\geq 0}$ is an ${\cal F}_i$-martingale.
For $q=0$, we recover the standard simple random walk.
We prove  the following.
\begin{theorem}
  \label{theo-main}
  For any $q\in (0,1)$, there exists $\sigma_+(q), \sigma_-(q)\in (0,1/2)$ and $c,C \in (0,\infty)$ such that for any $n$

  \begin{equation}
    \label{bounds}
     c n^{-\sigma_-(q)} < \sup_{u\in {\cal U}_q} P(S_n^u=0) < C n^{-\sigma_+(q)}
  \end{equation}
  and
\begin{equation}
\label{eq-final1}
\sigma_-(q)\to_{q\nearrow 1}0.
\end{equation}
\end{theorem}

Work related to ours (in the context of the control of diffusion processes)
has appeared in \cite{mc1}; more recently, the results in \cite{alexander} are related to the lower bound in Theorem \ref{theo-main}.
\section{Proofs}
Theorem \ref{theo-gen} (which immediately implies the upper bound in Theorem
\ref{theo-main})
%As is often the case in control problems, one direction (the upper bound)
is obtained by observing that any martingale  has to overcome
a (logarithmic number of) barriers in order to reach the target region,
and each such barrier can be overcome
only with (conditional on the history)
probability bounded away from $1$. The lower bound in Theorem \ref{theo-main},
on the
other hand, will be obtained by exhibiting an explicit control.
\subsection{Upper bound - Proof of Theorem \ref{theo-gen}}
Throughout this subsection, $\delta\in (0,1]$ is a fixed constant, and
$\{M_i\}_{i\geq 0}$ denotes a martingale adapted
to a filtration ${\cal F}_i$, satisfying the condition
\begin{equation}
  \label{eq-deltaLB}
  E\left( (M_{i+1}-M_i)^2|{\cal F}_i \right)\geq \delta\,.
\end{equation}
We begin with an elementary lemma.
\begin{lemma}
  \label{lem-0}
  Assume that $M_0=0$, that
  \eqref{eq-deltaLB} holds, and that
  for some $h\geq 1$,
  $|M_{i+1}-M_i|\leq h$ almost surely. Fix
  \begin{equation}
    \label{eq-RS0}
    \ell\geq 24 h^2/\delta\,.
  \end{equation}
  Let $\tau=\min\{i: |M_i|\geq h\}$.
Then,
\begin{equation}
  \label{eq-lem0}
  P(M_\tau \geq h, \tau\leq \ell)\geq \frac16\,.
\end{equation}
\end{lemma}
\textit{Proof of Lemma \ref{lem-0}.}
By \eqref{eq-deltaLB}, the process $\{M_i^2-\delta i\}$ is a sub-martingale,
hence
$$0\leq E(M^2_{\tau\wedge \ell}-\delta (\tau\wedge\ell))
\leq 4h^2-\delta E(\tau\wedge\ell)\,,$$
where the bound on the increments of $\{M_i\}$ was used in the last inequality.
It follows that $E(\tau\wedge\ell)\leq 4h^2/\delta$, and therefore,
\begin{equation}
  \label{eq-RS1}
  P(\tau\geq \ell)\leq 4h^2/\delta\ell\leq \frac16\,,
\end{equation}
where \eqref{eq-RS0} was used in the second inequality. On the other hand,
using again that increments of $\{M_i\}$ are bounded by $h$,
$$0=EM_\tau\leq 2h P(M_\tau\geq h)-h P(M_\tau\leq -h),$$
which implies that $P(M_\tau\leq -h)\leq 2/3$. Combining this with
\eqref{eq-RS1} yields the lemma.
\qed

We have the following corollary.
\begin{lemma}
  \label{lem-2}
  Let $H,L>0$ and let $K$ be a positive integer so that $H^2\leq \delta KL/24$.
  Assume \eqref{eq-deltaLB}, $M_0=0$, and that
  \begin{equation}
    \label{eq-longtube}
    |M_{i+1}-M_i|\leq \frac{H}{K}\,, \quad \mbox{almost surely}
  \end{equation}
  Let $\tau_H=\min\{i: M_i\geq H\}$.
  Then,
  \begin{equation}
    \label{eq-lem2}
    P(\tau_H\leq L)\geq \left(\frac16\right)^K\,.
  \end{equation}
\end{lemma}
\textit{Proof of Lemma \ref{lem-2}.}
Set $\ell=L/K$, $h=H/K$, and iterate Lemma \ref{lem-0} $K$ times.
\qed

Combining Lemma \ref{lem-0} and Lemma \ref{lem-2}, one obtains the following.
\begin{lemma}
  \label{lem-3}
  Let $L>0$ be a positive integer. Set $H=3\sqrt{L}$
  and let $K$ be a positive integer so that $K\delta\geq 216$.
  %$H^2\leq \delta KL/24$.
  Assume \eqref{eq-deltaLB}, $M_0=z$,  \eqref{eq-longtube} and
  \begin{equation}
    \label{eq-UBQV}
     E\left( (M_{i+1}-M_i)^2|{\cal F}_i \right)\leq 1\,.
   \end{equation}
  Let
$${\cal D}=\{ (x,i)\in \mathbb{R}\times \mathbb{Z}_+:
i\in [L,2L], |x|\leq H/3\}\,.$$
Then,
\begin{equation}
  \label{eq-D}
  P(\{(M_i,i)\}_{i=L}^{2L}\cap {\cal D}=\emptyset)\geq \frac12\cdot \left(
  \frac16\right)^K\,.
\end{equation}
\end{lemma}
\textit{Proof of Lemma \ref{lem-3}.}
It is enough to consider $z\geq 0$.
Let $\bar \tau_H=\min\{i: M_i\geq z+H\}$. Note that the condition on
$K$ ensured that $H^2\leq \delta K L/24$.
By Lemma \ref{lem-2},
$$P(\bar \tau_H\leq L)\geq (1/6)^K.$$ On the other hand, by Doob's inequality
and \eqref{eq-UBQV},
on the event $\bar \tau_H\leq L$,
$$P(\sup_{i\leq 2L}|M_{i+\bar \tau_H}-M_{\bar\tau_H}|\geq 2H/3|
{\cal F}_{\bar \tau_H})\leq \frac{2L}{(2H/3)^2}=\frac12\,.
$$
Combining the last two displays completes the proof.
\qed

We can now begin to construct the barriers alluded to above.
Fix $n>0$ and set $V_{m,n}=[m,n]\cap \ZZ$, $R_n=[-n,n]$.
Write $V_n=V_{0,n}$ and $B_{j,n}=V_{(1-2^{-j})n,n}$.
Define the following nested subsets of $\reals\times  V_n$:
$$ D_0=\reals\times [0,n],
D_{i}= R_{[2^{-i/2}\sqrt{n}]}\times B_{i,n}\,.
$$
Let $N_0=\max\{i: 2^{-i/2}\sqrt{n}\geq n^\beta\}$, i.e.
$$N_0\geq \frac{1/2-\beta}{\log \sqrt{2}}\log n-\frac{216/\delta}{\log \sqrt{2}}-1\,.$$

Let $\tau_0=0$ and for $i\geq 1$ set $\tau_i=\min\{t>\tau_{i-1}:
(M_{t},t)\in D_i\}$. A direct consequence of
Lemma \ref{lem-3} is the following.
\begin{lemma}
  \label{lem-tau}
  There exists a constant $c=c(\delta)>0$ so that
  on the event $\tau_i<n$, and with $i\leq N_0-1$, one has
  $$P(\tau_{i+1}\geq n|{\cal F}_{\tau_i})\geq c, a.s.
  \,.
$$
\end{lemma}
(The choice of $N_0$ ensured that in the applications of Lemma \ref{lem-3}
for any $i\leq N_0-1$, the condition
\eqref{eq-longtube} holds.)

\noindent
\textit{Proof of Theorem \ref{theo-gen}.}
It is clearly enough to consider $z=0$ with arbitrary $M_0$. Adjusting $C$ if necessary, we may and will assume that $N_0>1$.
Note that
$\{|M_n|\leq n^\beta\}\subset
\{\tau_{N_0}\leq n\}$ and therefore, by Lemma \ref{lem-tau},
$$ P(|M_n|\leq  n^\beta)\leq  (1-c)^{N_0-1}\,.$$
This yields the theorem.
\qed

\subsection{Proof of  Theorem \ref{theo-main}}
The upper bound in \eqref{bounds} is a consequence of Theorem \ref{theo-gen}.
We thus need only to consider the lower bound in  \eqref{bounds}, and
the claim \eqref{eq-final1}.

First note that the simple control $u_i=q$
%To avoid parity issues
already achieves
the lower bound with exponent $\sigma_-(q)=1/2$. Thus,
what we need to show is that for any $q>0$ there is a
(polynomially) better control and that as $q \to 1$ we can
achieve an exponent close to 0.
Toward this end,
we use two very simple controls, that
are not approximation of the optimal control. See Section \ref{sec-concl}
for further comments on this point.

We begin with the following a-priori estimate.
\begin{lemma}
  \label{lem-ori2}
For any $q>0$ there exist $\alpha>0$, $\beta>0$, $K_0>0$  and
$\eps>0$ such that for any $K>K_0$ there is a $q$-admissible
control such that
$$\sum_{x=-K}^K P_x (S_{\alpha K^2}^u =y) > 1+\eps \, ,$$
for any $y\in [-\beta K,\beta K]$.
\end{lemma}
\textit{Proof of Lemma \ref{lem-ori2}:}
The control we take is slow inside $[-\beta K, \beta K]$ and fast outside, i.e. we take $u_i=q$ for $|S_i^u| \le \beta K$ and $u_i=0$ for
$|S_i^u| > \beta K$. We claim that given any $q>0$,
using this control with $\alpha>0$ and $\beta>0$ small enough and $K>K_0$ with $K_0$ large enough
will satisfy the conclusion of the lemma with some $\eps>0$.

Our control does not change with time, it is a reversible
Markov chain with weights $w_{x,x+1}=1$ and
$w_{x,x}=0$ for $|x| > \beta K$ and
$w_{x,x}=2q/(1-q)$ for $|x|\le \beta K$. Its reversing
measure is thus $\pi_x=2$ for $|x| > \beta K$ and
$\pi_x=2/(1-q)$ for $|x|\le \beta K$.

Using reversibility we get
$$P_x(S_{\alpha K^2}^u =y)= P_y(S_{\alpha K^2}^u =x) \frac{\pi_y}{\pi_x}$$

Thus,

$$\sum_{x=-K}^K P_x (S_{\alpha K^2}^u =y)=
\sum_{x=-K}^K P_y(S_{\alpha K^2}^u =x) \frac{\pi_y}{\pi_x}$$

$$=\frac{1}{1-q} \Big[\sum_{x=-K}^{-\beta K} P_y(S_{\alpha K^2}^u =x) +\sum_{x=\beta K}^{K} P_y(S_{\alpha K^2}^u =x)\Big] + \sum_{x=-\beta K}^{\beta K} P_y(S_{\alpha K^2}^u =x)$$

$$=\frac{1}{1-q} \Big[P_y(S_{\alpha K^2}^u \in [-K,K])- q P_y(S_{\alpha K^2}^u \in [-\beta K, \beta K])\Big]$$

Now, the probability that a simple random walk
will get to a distance of more than $K/2$ in $\alpha K^2$ steps
tends to 0 as $\alpha$ tends to 0, uniformly in $K$. Obviously,
this applies also for our controlled random walk (which is sometimes lazy),
hence by choosing small enough $\alpha$ we can
guarantee that for any $K>0$
and
any $y\in [-K/2,K/2]$ we have $P_y(S_{\alpha K^2}^u \in [-K,K])>1-q$.

Having fixed $\alpha$, we now claim that
%as $\beta\to 0$ we have
\begin{equation}
\label{eq-betaK}
 \limsup_{K_0\to \infty} \limsup_{\beta\to 0}
\sup_{K>K_0}
\sup_{y\in [-\beta K,\beta K]}
P_y(S_{\alpha K^2}^u \in [-\beta K, \beta K])= 0\,.
\end{equation}
%uniformly in $K$.
Indeed, by \cite[Corollary 14.5]{woess}, there exists a constant
$C(q)$ so that
%for any infinite countable state space
%reversible Markov chain we have
%$$p^t(x,y)\le \frac{\pi_x}{\pi_* \sqrt{t}} \, ,$$
$$p^t(x,y)\le \frac{C(q)}{ \sqrt{t}} \, ,$$
for any two states $x$ and $y$. (The bound in \cite{woess} is valid for any
random walk on an infinite graph with bounded degree and bounded above and below
conductances, see \cite[Pg. 40]{woess}; Note that while the bound is stated for
discrete time chains, it can also be transferred without much effort to
the continuous
time setting. See e.g. \cite[Theorem 2.14 and Proposition 3.13]{Kumagai}.)
%\textbf{[YUVAL: add reference, fix inequality?  what is $\pi_*$?]}
%there's probably some constant there]}

Plugging $t=\alpha K^2$, we get
$$P_y(S_{\alpha K^2}^u \in [-\beta K, \beta K] <
\frac{C(q)(2\beta K +1)}{\sqrt{\alpha}K} \, ,$$
%\frac{2\beta K +1}{(1-q)\sqrt{\alpha}K} \, ,$$
which tends to 0 when $\beta\to 0$ and $K\to\infty$ in the order
prescribed in \eqref{eq-betaK}.
%$\beta\to 0$ and $K\to \infty$.

Thus, by choosing small enough $\beta$ and large enough
$K_0$ we can have
$$P_y(S_{\alpha K^2}^u \in
[-K,K])- q P_y(S_{\alpha K^2}^u \in [-\beta K, \beta K]) > 1-q $$
uniformly for all  $K>K_0$ and we are done.
\qed

\noindent
\textit{Proof of the lower bound in Theorem \ref{theo-main}:}
Fix $q>0$ and choose $\alpha, \beta, K_0$ and
$\eps$ according to Lemma \ref{lem-ori2}.

Let $L=\lfloor - \log (T/K_0^2) / 2 \log \beta \rfloor$ and let
$T_\ell=T-\alpha K_0^2 \sum_{i=1}^\ell \beta^{-2\ell}$, for
$\ell=1,\ldots,L$ and $T_0=T$. For time $t=0,\ldots,T_L$ we
use the control $u_t=q$.
Notice that $T_L \approx T$ so standard estimates for
lazy random walk show that there exists a constant $c>0$,
independent of $T$, such that $P_0(S_{T_L}^u=y)>c T^{-1/2}$,
for any $y\in [-K_0 \beta^{-L},K_0 \beta^{-L}]\subset [-T^{1/2},T^{1/2}]$.

For any $\ell=L,\ldots,1$,
from time $T_\ell$ to $T_{\ell-1}$ we use the
strategy provided by Lemma \ref{lem-ori2} for $K=K_0 \beta^{-\ell}$.
Applying Lemma \ref{lem-ori2} repeatedly, we see that for any
$\ell=L-1,\ldots,0$,
at time $T_\ell$ we have $P_0(S_{T_\ell}^u=y)>c (1+\eps)^{L-\ell}
T^{-1/2}$ for any $y\in [-K_0 \beta^{-\ell},K_0 \beta^{-\ell}]$.
In particular, we have
$$P_0(S_T^u=0)>c' (1+\eps)^L T^{-\frac12}=c (1+\eps)^{\frac{\log K_0}{\log \beta}} T^{-\frac12 -\frac{\log(1+\eps)}{2 \log \beta}} \, ,$$
showing that $\sigma_-(q)<1/2$.
This completes the proof of \eqref{bounds}. \qed

In preparation for the proof of \eqref{eq-final1}, we
provide an auxilliary estimate.
\begin{lemma}
  \label{lem-ori}
For any $\eps>0$ there exist $A$ and $q<1$
such that for any $K$ there is a $q$-admissible
control with the property that for any $x\in [-2K,2K]$ we have
$$P_x (S_{A K^2}^u \in [-K, K]) > 1 - \eps \, .$$
\end{lemma}
\textit{Proof of Lemma \ref{lem-ori}:}
%The control which is slow in $[-K,K]$ and fast outside
%should do the trick for any big enough $\alpha$ and $q$,
%but it is even easier to prove for the following control:
Let $A$ be so that for a simple random walk
on $\ZZ$ we have for any $K$,
\begin{equation} \label{ori1}
P_0(\tau_{2K}>A K^2)< \eps /2 ,
\end{equation}
where $\tau_{2K}$ is the first hitting time of $2K$.

Having chosen $A$, let $q<1$ be so big such that for a $q$-lazy
random walk (that is, a random walk with control $u_i=q$)
we have for any $K$,
\begin{equation} \label{ori2}
P_0(\tau_{\{-K,K\}}<A K^2)<\eps/2 ,
\end{equation}
where $\tau_{\{-K,K\}}$ is the first time of hitting either $K$ or $-K$.

We now define the control to be fast until the walk
hits $0$ and slow afterwards, i.e. we take $u_i=0$ for
$i<\tau_0:=\min\{n: S_n^u=0\}$ and $u_i=q$ for
$i\geq \tau_0$.
If the starting location $S_0^u$ is in $[-2K,2K]$,
then by (\ref{ori1}) with probability at least $1-\eps/2$ we hit
$0$ before time $A K^2$.
If that happens, then by (\ref{ori2}) with probability at least
$1-\eps/2$, the walk stays inside $[-K,K]$ until time $A K^2$. \qed

We can now complete the proof of Theorem \ref{theo-main}.\\
\textit{Proof of \eqref{eq-final1}:}
%$\lim_{q\nearrow 1}\sigma_-(q) =0$.}\\
Fix $\eps>0$ and choose $q$ and $A$ according to Lemma \ref{lem-ori}.

Let $L=\lfloor \log_4 (n/A)\rfloor$ and let
$T_\ell=T-A \sum_{i=0}^\ell 4^k$, for $\ell=0,\ldots,L$.
For time $0$ to $T_L$, we have
$$P_0(S_{T_L}^u\in [-2^L,2^L])>c$$
for some fixed $c>0$, regardless of the control.

For any $\ell=L,\ldots,1$, from time $T_\ell$ to $T_{\ell-1}$ we use the
strategy provided by Lemma \ref{lem-ori}
for $K=2^\ell$. Then with probability at least
$c (1-\eps)^L \approx n^{\log_4 (1-\eps)}$ we have $S_T=0$.
This yields the required lower bound.
\qed

\section{Concluding remark}
\label{sec-concl}
Motivated by the structure of the optimal control in the
continuous time-and-space analogue of our control problem, see
\cite{Atrok}, one could
attempt to improve on the lower bound in \eqref{bounds}
by using a bang-bang control of
the type $u_i=q$ if $(S_i^u,i)\in D\subset \ZZ\times [0,n]$ and
$u_i=0$ otherwise, where $D$ is a domain whose boundary is
determined by an appropriate (roughly
parabolic) curve. The analysis of that control
is somewhat tedious, and proceeding in that direction we have
only been able to show the lower bound in \eqref{bounds} with
%while we have been able to prove \eqref{bounds}
%using such analysis,
%we have not been able to refine it so as to show that
$\sigma_-(q)<1/2$ when $q$ is sufficiently large.
%close to $1$.
It would be interesting to check whether an analysis of
the dynamic programming
equation associated with the control problem, in line with its continuous
time analogue in \cite{Atrok,BaS}, could yield that estimate, and
more ambitiously, show the equality of $\sigma_-(q)$ and $\sigma_+(q)$
in \eqref{bounds}. 

One could also consider the dual problem of \textit{minimizing} the probability
of hitting $0$ at time $n$, that is, in the setup of Theorem \ref{theo-main}, 
of evaluating
\begin{equation}
    \label{bounds-min}
%     c n^{-\sigma_-(q)} < 
\inf_{u\in {\cal U}_q} P(S_n^u=0) \,.
\end{equation}
One can adapt the proof of the lower bound in Theorem \ref{theo-main}
(replacing in the sub-optimal control ``fast'' by ``slow'') to obtain
a polynomial upper bound in \eqref{bounds-min} that has exponent larger than 
$1/2$. Similarly (using the invariance principle for martingales),
one shows that there is $\alpha=\alpha(q)$ such tha
the controlled walk with $|S_0^u|<2K$ satisfies 
$|S_{\alpha K^2}|\leq K$ with positive (depending only on $q$ and
independent of $K$) probability, and from this a polynomal lower bound
in \eqref{bounds-min} follows. We omit further details.
%\textit{must}
%  \end{equation}
\\[2em]

\noindent
{\bf Acknowledgement} We thank Bruno Schapira for pointing out \cite{TDLR} to us, and Charlie Smart for his crucial comments \cite{Smart}.

\end{document}